\documentclass[a4paper]{article}

\usepackage[english]{babel}
\usepackage[utf8]{inputenc}
\usepackage{amsmath}
\usepackage{graphicx}
\usepackage[colorinlistoftodos]{todonotes}
\usepackage[affil-it]{authblk}

\usepackage{amssymb} 
\usepackage{epstopdf} 
\usepackage{amsmath} 
\usepackage{amsthm} 
\usepackage{mathrsfs}

\newtheorem{theorem}{Theorem}[section] 
\newtheorem*{mn_thm_1}{Theorem A}
\newtheorem*{mn_thm_2}{Theorem B}
\newtheorem*{mn_thm_3}{Theorem C}
\newtheorem{exmp}{Example}[section]

\newtheorem*{theorem-A}{Theorem A}

\newtheorem{remark}[theorem]{Remark}

\title{Maximal abelian normal subgroups in Nilpotent groups of class $2$}

\author{Satvik Goswami%
  \thanks{Undergraduate student, Electronic address \texttt{satvik.goswami@gmail.com}}}
\affil{Indian Institute of Science Education and Research, \\ 
Bhopal, India}

\author{Ashish Gupta%
  \thanks{Corresponding author, Electronic address: \texttt{a0gupt@gmail.com}}}
\affil{Department of Mathematics and  
Statistics,\\
University of Melbourne, \\
Parkville, Melbourne, \\
Australia.}

 \affil{
Department of Mathematics,\\ Indian Institute of Science Education and Research, \\Bhopal, India.}

\date{}

\begin{document}

\maketitle

\begin{abstract}
A maximal abelian normal subgroup $A$ in a  nilpotent group $N$ is self-centralizing. This makes their role an important one in determining the structure of the nilpotent group. For example if $A$ is finite then $N$ is also finite.

In the free nilpotent group of class $2$ such a group is always a cyclic extension of the center. However this need not be true in an arbitrary finitely generated class two nilpotent group. We show that three generator nilpotent groups violating this condition can be easily found. However in the case when there are four or more generators certain conditions must hold.     
\end{abstract}

\section{Introduction}
Nilpotent groups play an important role in the classification of groups and that of Lie groups. 
Whereas the simplest nilpotent groups are the abelian groups, the nilpotent groups which are at the next level of complexity are the class two nilpotent groups. The equation $[[x,y],z] = 1$ is satisfied identically in such a group where $[x, y] = x^{-1}y^{-1}xy$. For example, groups of $3 \times 3$  unitriangular matrices over a commutative ring $R$ are of this type.\\
A maximal normal abelian subgroup $A$ of a nilpotent group $N$ is known to exert a strong influence on $N$. For example if $A$ is finite then $N$ is also finite. This is due to the self-centralising property of $A$. It is easy to see that such a subgroup will always contain the center. 

Consider the free nilpotent group $F_n$ on $n$ generators $x_1, \cdots, x_n$. In such group the center is freely generated by the $\frac{n(n - 1)}{2}$ commutators $[x_i, x_j]$ and $F_n$ modulo its center is free abelian with basis given by the images of the generators $x_i$. In such a group a maximal abelian normal subgroup is an extension of the center by a cyclic subgroup. This property may not be preserved in passing to quotients of $F_n$ by subgroups generated by central elements. Out aim in this paper is to investigate this question. More precisely, consider the case $n = 3$ and the group 
\[ G(a_1, a_2, a_3) :=  H_3/\langle [x_1, x_2]^a_1[x_2, x_3]^a_2[x_1,x_3]^a_3\rangle \] 
so that $(a_1, a_2, a_3) \ne (0,0,0)$. As the following result shows in this case maximal abelian normal subgroups in $G(a_1, a_2, a_3)$ have (torsion-free) rank at least two over the center.    

\begin{mn_thm_1}
Let $H_3$ be a free nilpotent group of rank 3 and class 2 with basis $x_i$, $1 \leq i \leq 3$. Let \[ H = F_3/\langle [x_1,x_2]^{\ a_1}[x_2,x_3]^{ a_2}[x_1,x_3]^{a_3} \rangle, \] 
where $a_i \in \mathbb{Z} \setminus \{0\}$. Then for a maximal normal abelian subgroup $A$ in $H$,
\[ rk(A) \ge rk(\zeta H) + 2. \] 
\end{mn_thm_1}
For $n \ge 4$ we find that $F_n/C$ can have quotients by central elements in which maximal abelian normal subgroups are cyclic extensions of the center (and thus have rank one over the center).
In the case of four generators the following holds.

\begin{mn_thm_2}
Let $C < F_4$ be the cyclic subgroup $\langle \prod_{1\leq i < j\leq 3} \left[x_i,x_j\right]^{a_{ij}} \rangle$ with $a_{ij}\in \mathbb{Z} \setminus \{0\}$ and
\[ \epsilon := \frac{|\prod a_{ij}|}{\prod a_{ij}} \]
Then for $F_4/C$ to have abelian subgroups containing the center and having rank greater than one over the center, the following condition must hold
\begin{align*}
\epsilon\left(\frac{1}{a_{14}a_{23}} + \frac{1}{a_{12}a_{34}}\right) > \epsilon\left(\frac{1}{a_{13}a_{24}}\right).
\end{align*}
\end{mn_thm_2}

\begin{mn_thm_3}
Let $F_n$ denote the free nilpotent group of class $2$ and rank $n \geq 3$ with basis $\{x_i, 1 \leq i \leq n\}$. Let $C < \zeta F_n$ be the cyclic subgroup 
\[ C : =  \langle \prod_{1\leq i < j\leq n} \left[x_i,x_j\right]^{a_{ij}} 
\rangle \] where 
$a_{ij} \in \mathbb{Z} \setminus \{0\}$ and
\[ \epsilon := \frac{|\prod a_{ij}|}{\prod a_{ij}}. \]
Then for $F_n/C$ to have abelian subgroups containing the centre and having rank greater than one over the centre, the following conditions must hold
\begin{align*}
\epsilon\left(\frac{1}{a_{k_1k_4}a_{k_2k_3}} + \frac{1}{a_{k_1k_2}a_{k_3k_4}}\right) > \epsilon\left(\frac{1}{a_{k_1k_3}a_{k_2k_4}}\right).
\end{align*}
for each $\{k_1,k_2,k_3,k_4\} \subseteq \{1,...,n\}$ such that $1<k_1<k_2<k_3<k_4<n.$

\end{mn_thm_3}

\section{The three generator case}

We shall prove the following theorem.

\begin{mn_thm_1} \label{Thm_A}
Let $F_3$ be a free nilpotent group of rank 3 and class 2 with basis $x_i$, $1 \leq i \leq 3$. Let \[ H = F_3/\langle [x_1,x_2]^{\ a_1}[x_2,x_3]^{ a_2}[x_1,x_3]^{a_3} \rangle, \] 
where $a_i \in \mathbb{Z} \setminus \{0\}$. Then for a maximal normal abelian subgroup $A$ in $H$,
\[ rk(A) \ge rk(\zeta H) + 2. \] 
\end{mn_thm_1}

\begin{proof}
Let 
\[ C : = \langle [x_1,x_2]^{\ a_1}[x_2,x_3]^{ a_2}[x_1,x_3]^{a_3} \rangle.\]
We shall prove by constructing a subgroup $A$ of the required type. For $1 \leq i  \leq 2,$ let $\alpha_i \in F_3$ be the element with an expression: $\alpha_i = \prod_{j = 1}^3 x_j^{m_{ij}}$. Then 
\[ [\alpha_1, \alpha_2] = \prod_{i < j} [x_i, x_j]^{d_{ij}}, \]
where, \[  d_{ij}  = 
\det\begin{bmatrix}
m_{1k}& m_{2k}\\
m_{1l}& m_{2l}
\end{bmatrix}. \]
The above equations holds because the commutators $[x_i, x_j]$ are central in $H_3$ (see \cite[Section 5.1.5]{Ro}).
In what follows round parentheses will denote the greatest common divisor. Let $w_1$, $w_2 \in \mathbb{Z}$ be a solution of:
\begin{align} 
\frac{a_1(a_2,a_3)}{(a_1,a_2,a_3)}X+ \frac{a_2(a_1,a_3)}{(a_1,a_2,a_3)}Y- \frac{a_3(a_1,a_2)}{(a_1,a_2,a_3)} &= 0 \nonumber
\end{align}
Then it is easily checked that the assignments 
\begin{align*}
m_{11} &:= (a_1,a_3)w_2, \ \ \ \ \ 
m_{12} := (a_1,a_2), \\ 
m_{13} &:= (a_2,a_3)w_1, \ \ \ \ \ 
m_{21} :=-a_1/(a_1,a_2), \\  
m_{22} &:= 0, \ \ \ \ \  \ \ \ \ \ \ \ \ \ \ \ 
m_{23} :=a_2/(a_1,a_2)
\end{align*}
yield
\[ \det\begin{bmatrix}
m_{11}& m_{21}\\
m_{12}& m_{22}
\end{bmatrix}
= a_1, \  \det\begin{bmatrix}
m_{12}& m_{22}\\
m_{13}& m_{23}
\end{bmatrix}= a_2 \ \ \mbox{and} \ \  \det\begin{bmatrix}
m_{11}& m_{21}\\
m_{13}& m_{23}
\end{bmatrix}
= a_3. \] 
By our construction $[\alpha_1,\alpha_2] \in C.$ Hence $ [\alpha_{1}C,\alpha_{2}C] = 1$ in $F_3/C$, having generators ${x_1}C,\ {x_2}C$ and ${x_3}C$. Now let $\mathcal Z$ denote the center of ($F_3/C)$ and define
\[A := \langle {\alpha_1}C,\ {\alpha_2}C,\ \mathcal Z ({F_3}/C) \rangle \]
We claim that $A/ \mathcal Z$ is freely generated by ${\alpha_1}C$ and ${\alpha_2}C$ and so is free abelian of rank $2$. To this end let 
$(m_1,m_2) \in \mathbb Z^2$ so that  
\[ [({\alpha_1}C)\mathcal {Z}]^{m_1} [({\alpha_2}C)\mathcal {Z}]^{m_2} = 1 \]
which implies that $({\alpha_1}C)^{m_1}({\alpha_2}C)^{m_2} \in \mathcal {Z}$, or in other words
\begin{align}
{\alpha_1}^{m_1}{\alpha_2}^{m_2}C \in \mathcal Z.
\end{align}
We wish to show that $m_1 = m_2 = 0.$ The last equation above means that each of ${x_1}C,\ {x_2}C$ and ${x_3}C$ commutes with ${\alpha_1}^{m_1}{\alpha_2}^{m_2}C.$ In particular, we have
\[ [{\alpha_1}^{m_1}{\alpha_2}^{m_2}C,\ {x_1}C] = 1. \]
Since $F_3$ is nilpotent of class $2$ we can write above equation as
\begin{align}
[\alpha_{1}C, x_{1}C]^{m_1}[\alpha_{2}C, x_{1}C]^{m_2} &= 1. 
\end{align}
Substituting ${\alpha_1} := {x_1}^{m_{11}}{x_2}^{m_{12}}{x_3}^{m_{13}}$ and ${\alpha_2} := {x_1}^{m_{21}}{x_2}^{m_{22}}{x_3}^{m_{23}}$ in the last equation we get
\[ [{x_1}^{m_{11}}{x_2}^{m_{12}}{x_3}^{m_{13}}, x_1]^{m_1}[{x_1}^{m_{21}}{x_2}^{m_{22}}{x_3}^{m_{23}}, x_1]^{m_2} \in C \]
On simplification using \cite[Section 5.1.5]{Ro} we easily get
\begin{align}\label{simpf_comtr_reln}
[{x_1},{x_2}]^{-(m_{12}m_{1}+m_{22}m_{2})}[{x_1},{x_3}]^{-(m_{13}m_{1}+m_{23}m_{2})} \in C
\end{align}
Comparison of equation $(\ref{simpf_comtr_reln})$ with $C = \langle [x_1,x_2]^{a_1}[x_2,x_3]^{ a_2}[x_1,x_3]^{a_3} \rangle$ show that $a_2 = 0$. But this contradicts the assumption on $a_2$.

\end{proof}
\begin{remark}
We contend that for quotients $F_3/\langle \left [x_1, x_2\right]^{\ a_1}\left[x_2,x_3\right]^{\ a_2}\left[x_1,x_3\right]^{\ a_3}\rangle $ a maximal normal abelian subgroup has rank at most $2$ over the center so that the inequality in Theorem A is an equality.
\end{remark}

\section{The four generator case}

\begin{mn_thm_2}\label{Thm_B}
Let $C < F_4$ be the cyclic subgroup $\langle \prod_{1\leq i < j\leq 4} \left[x_i,x_j\right]^{a_{ij}} \rangle$ with $a_{ij}\in \mathbb{Z} \setminus \{0\}$ and
\[ \epsilon := \frac{|\prod a_{ij}|}{\prod a_{ij}} \]
Then for $F_4/C$ to have abelian subgroups containing the center and having rank greater than one over the center, the following condition must hold
\begin{align*}
\epsilon\left(\frac{1}{a_{14}a_{23}} + \frac{1}{a_{12}a_{34}}\right) > \epsilon\left(\frac{1}{a_{13}a_{24}}\right).
\end{align*}
\end{mn_thm_2}
\begin{proof}
Let
\[ \alpha_i=\prod_{1\leq j\leq 4} x_j^{m_{ij}} \in F_4,\]
where $1\leq i\leq 2$ with $[\alpha_1,\alpha_2]\in C.$ We have  
\[ [\alpha_1,\alpha_2]=\prod_{1\leq i < j\leq 4} \left[x_i,x_j\right]^{d_{ij}}, \]
where \[ d_{ij} =
det\begin{bmatrix}
m_{1i}& m_{2i}\\
m_{1j}& m_{2j}
\end{bmatrix}. \]  Hence if $[\alpha_1,\alpha_2]\in C$ then $d_{ij}= la_{ij}$ for some $l\in \mathbb{Z}$.\\
Writing this explicitly we have,
\begin{align}
\label{a_12}
m_{11}m_{22}-m_{12}m_{21} &=la_{12},\\
\label{a_23}
m_{12}m_{23}-m_{13}m_{22} &=la_{23},\\ 
\label{a_13}
m_{11}m_{23}-m_{13}m_{21} &=la_{13},\\
\label{a_14}
m_{11}m_{24}-m_{21}m_{14} &=la_{14},\\ 
\label{a_24}
m_{12}m_{24}-m_{22}m_{14} &=la_{24},\\
\label{a_34}
m_{13}m_{24}-m_{23}m_{14} &=la_{34}
\end{align}
By equations ($\ref{a_12})$ and ($\ref{a_23}$), we get
\begin{align}
\label{m_11_m_22}
m_{11}m_{22} &=la_{12}+m_{12}m_{21},\\
\label{m_13_m22}
m_{13}m_{22} &=m_{12}m_{23}-la_{23}
\end{align}
Multiplying equation ($\ref{a_13})$ by $m_{22},$ we have
\begin{align}\label{a_13m_22}
(m_{11}m_{22})m_{23}-(m_{13}m_{22})m_{21} &=la_{13}m_{22}
\end{align}
Combining (\ref{m_11_m_22}) and (\ref{m_13_m22}) with (\ref{a_13m_22}) we obtain 
\begin{align}\label{dio_fourg_1}
a_{12}m_{23}+a_{23}m_{21} &=a_{13}m_{22}
\end{align}
In a similar way multiplying equation ($\ref{a_13})$ by $m_{12},$ we have
\begin{align}\label{a_13m_12}
(m_{12}m_{23})m_{11}-(m_{12}m_{21})m_{13} &=la_{13}m_{12}
\end{align}
Now substituting ($m_{12}m_{21})$ and ($m_{12}m_{23})$ from equations ($\ref{a_12})$ and ($\ref{a_23})$ respectively in ($\ref{a_13m_12})$ and cancelling off $l$, we obtain
\begin{align}\label{dio_fourg_2}
a_{12}m_{13}+a_{23}m_{11} &=a_{13}m_{12}
\end{align}
Equations ($\ref{dio_fourg_1})$ and ($\ref{dio_fourg_2})$ are regarded as diophantine equations in variables $m_{i1}, m_{i2}$ and $m_{i3}$ have the same form, namely,\\
\begin{align}\label{fourg_a_13_m_i2}
a_{12}m_{i3}+a_{23}m_{i1}=a_{13}m_{i2},\ \  1\leq i\leq 2
\end{align}
Similarly considering other combinations of equations (\ref{a_12}) -- (\ref{a_34}) we obtain six more diophantine equations:
\begin{align*}
a_{34}m_{i2}+a_{23}m_{i4}=a_{24}m_{i3},\\
a_{34}m_{i1}+a_{13}m_{i4}=a_{14}m_{i3},\\
a_{24}m_{i1}+a_{12}m_{i4}=a_{14}m_{i2},
\end{align*}
where $1\leq i\leq 2.$ Writing these equations explicitly we have
\begin{align*}
a_{23}m_{11}+a_{12}m_{12}-a_{13}m_{13}=0\\
a_{24}m_{11}+a_{12}m_{14}-a_{14}m_{12}=0\\
a_{34}m_{11}+a_{13}m_{14}-a_{14}m_{13}=0\\
a_{34}m_{12}+a_{23}m_{14}-a_{24}m_{13}=0\\
a_{23}m_{21}+a_{12}m_{23}-a_{13}m_{22}=0\\
a_{24}m_{21}+a_{12}m_{24}-a_{14}m_{22}=0\\
a_{34}m_{21}+a_{13}m_{24}-a_{14}m_{23}=0\\
a_{34}m_{22}+a_{23}m_{24}-a_{24}m_{23}=0
\end{align*}
But this means that $\det(M) = 0$ where 
\[ M := 
\begin{pmatrix}
A & 0\\
0 & A
\end{pmatrix} \] 
and 
\[ A := 
\begin{pmatrix}
a_{23} & -a_{13} & a_{12} & 0\\
a_{24} & -a_{14} & 0 & a_{12}\\
a_{34} & 0 & -a_{14} & a_{13}\\
0 & a_{34} & -a_{24} & a_{23}
\end{pmatrix} \] 
Since $\det(M) = 0$ hence $\det(A) = 0$ and we thus get
\[ (a_{12}^2a_{34}^2 + a_{13}^2a_{24}^2 + a_{14}^2a_{23}^2) + 2(a_{12}a_{14}a_{23}a_{34} - a_{13}a_{14}a_{23}a_{24} - a_{12}a_{13}a_{24}a_{34}) = 0 \]
Since 
\[(a_{12}^2a_{34}^2 + a_{13}^2a_{24}^2 + a_{14}^2a_{23}^2) > 0,\ \ \forall \ a_{ij}\in \mathbb{Z} \setminus \{0\}\]
This implies
\[ a_{12}a_{14}a_{23}a_{34} - a_{13}a_{14}a_{23}a_{24} - a_{12}a_{13}a_{24}a_{34} < 0\]
\begin{align}\label{cond_pre2}
a_{13}a_{14}a_{23}a_{24} + a_{12}a_{13}a_{24}a_{34} > a_{12}a_{14}a_{23}a_{34}  
\end{align}
On multiplying (\ref{cond_pre2}) by $\frac{|\prod a_{ij}|}{\left(\prod a_{ij}\right)^2}$ and denoting 
\[ \epsilon := \frac{|\prod a_{ij}|}{\prod a_{ij}} \]
we get the condition of the theorem:
\begin{align*}
\epsilon\left(\frac{1}{a_{14}a_{23}} + \frac{1}{a_{12}a_{34}}\right) > \epsilon\left(\frac{1}{a_{13}a_{24}}\right).
\end{align*}
\end{proof}
\begin{exmp}
For $F_4/C$ where
\[ C := \left [x_1, x_2\right]\left [x_1, x_3\right]\left [x_1, x_4\right]\left [x_2, x_3\right]^{-1}\left [x_2, x_4\right]\left [x_3, x_4\right]^{-1} \] 
it can be easily checked the above condition is violated. Hence a maximal abelian normal subgroup in this group is always a cyclic extension of the center.

\end{exmp}
\section{The $n$ generator case}

We have the following result  concerning the situation for central quotients of the free nilpotent group on $n$ generators.

\begin{mn_thm_3}\label{Thm_C}
Let $F_n$ denote the free nilpotent group of class $2$ and rank $n \geq 3$ with basis $\{x_i, 1 \leq i \leq n\}$. Let $C < \zeta F_n$ be the cyclic subgroup 
\[ C : =  \langle \prod_{1\leq i < j\leq n} \left[x_i,x_j\right]^{a_{ij}} 
\rangle \] where 
$a_{ij} \in \mathbb{Z} \setminus \{0\}$ and
\[ \epsilon := \frac{|\prod a_{ij}|}{\prod a_{ij}}. \]
Then for $F_n/C$ to have abelian subgroups containing the centre and having rank greater than one over the centre, the following conditions must hold
\begin{align*}
\epsilon\left(\frac{1}{a_{k_1k_4}a_{k_2k_3}} + \frac{1}{a_{k_1k_2}a_{k_3k_4}}\right) > \epsilon\left(\frac{1}{a_{k_1k_3}a_{k_2k_4}}\right)
\end{align*}
for each $\{k_1,k_2,k_3,k_4\} \subseteq \{1,...,n\}$ such that $1<k_1<k_2<k_3<k_4<n.$
\end{mn_thm_3}

\begin{proof}
Let
\[ \alpha_i=\prod_{1\leq j\leq n} x_j^{m_{ij}}\ \in F_n,\ \ \ \ 1\leq i\leq 2 \]
with $[\alpha_1,\alpha_2]\in C$. We have 
\[ [\alpha_1,\alpha_2]=\prod_{1\leq i < j\leq n} \left[x_{i},x_{j}\right]^{d_{ij}}\]
where 
\[ d_{ij}=det\begin{bmatrix}
m_{1i}& m_{2i}\\
m_{1j}& m_{2j}
\end{bmatrix}. \] 
Hence if $[\alpha_1,\alpha_2]\in C$, then \begin{equation} \label{eqn_for_C}
d_{ij}=la_{ij} \end{equation}
for some $l\in \mathbb{Z}$. In particular, the equations (\ref{eqn_for_C}) hold for every subset $\{k_1, k_2, k_3, k_4 \}$ of $\{ 1, \cdots, n \}$. 
So on comparing with the four generator case (Theorem B) we  easily arrive at the assertion of the theorem.

\end{proof}

\section*{Acknowledgements}
The first theorem was completed when the second author was a masters student at the University of Melbourne. This author gratefully acknowledges the unfailing support of his supervisor Dr. John R. J. Groves and also the university for financial support.

\end{document}